\documentclass[aap,MSNbibl,dvips]{arximspdf}
\usepackage{graphics}

\doi{10.1214/09-AAP629}
\volume{20}
\issue{2}
\pubyear{2010}
\firstpage{640}
\lastpage{659}

\makeatletter
\newtheorem{Lem}[Df]{Lemma}
\newtheorem{Th}[Df]{Theorem}

\newproclaim{Df}{Definition}[section]
\newproclaim{Ex}[Df]{Example}
\newproclaim{Com}[Df]{Comment}

\DeclareMathAlphabet\mathcaligr{OMS}{cmsy}{m}{n}
\renewcommand{\mathcal}{\mathcaligr}
\makeatother

\begin{document}
\begin{frontmatter}

\title{Optimal detection of a change-set in\\ a spatial Poisson process}
\runtitle{Detection of a change-set}

\begin{aug}
\author[a]{\fnms{B. Gail} \snm{Ivanoff}\thanksref{t1}\ead[label=e1]{givanoff@uottawa.ca}}\and
\author[b]{\fnms{Ely} \snm{Merzbach}\corref{}\ead[label=e2]{merzbach@macs.biu.ac.il}}
\thankstext{t1}{Supported by a grant from the Natural Sciences and
Engineering Research Council of Canada.}
\runauthor{B. G. Ivanoff and E. Merzbach}
\affiliation{University of Ottawa and Bar-Ilan University}
\address[a]{Department of Mathematics \& Statistics\\
University of Ottawa\\
585 King Edward\\
Ottawa, Ontario K1N 6N5\\
Canada\\
\printead{e1}}
\address[b]{Department of Mathematics\\
Bar-Ilan University\\
52900 Ramat-Gan\\
Israel\\
\printead{e2}}
\end{aug}

\received{\smonth{3} \syear{2008}}
\revised{\smonth{7} \syear{2009}}

%
\begin{abstract}
We generalize the classic change-point problem to a
``change-set'' framework: a spatial Poisson process changes
its intensity on an unobservable random set. Optimal
detection of the set is defined by maximizing the expected
value
of a gain function. In the case that the unknown
change-set is defined by a locally finite set of incomparable
points, we present a sufficient condition for optimal
detection of the set using multiparameter martingale techniques. Two
examples are discussed.
\end{abstract}



\end{frontmatter}

\section{Introduction}\label{Introduction}

In this paper, we consider the multiparameter version of the classic
optimal detection problem; the goal is to detect the
occurrence of a random set on which an observable Poisson
process changes its intensity. To be precise, we let $N=\{N_t,t\in
\mathbf{R}^2_+\}$ be a nonexplosive point process defined on the
positive quadrant of the plane and let
$\{\tau_n\}$ be its jump points, numbered in some arbitrary way. Then
$N_t=\sum_{n=1}^\infty I_{\{\tau_n\leq t\}}$ (cf. \cite{IM1}).
Here, ``$\leq$'' denotes the usual partial order on $\mathbf{R}_+^2{}\dvtx{} s=(s_1,s_2)\leq t=(t_1,t_2)\Leftrightarrow s_1\leq
t_1, s_2\leq t_2$. On some random set $\xi$, the intensity of
$N$ changes from $\mu_0$ to $\mu_1$, where $0<\mu_0<\mu_1$:
specifically, given
$\xi$, $N$ is a Poisson process with intensity
\[
\mu_0I_{\{t\notin\xi\}}+\mu_1I_{\{t\in\xi\}}
=\mu_0+(\mu_1-\mu_0)I_{\{t\in\xi\}}.
\]
The problem is that the ``change-set'' $\xi$ is
unobservable and we must detect $\xi$ as well as
possible, given our observation of the point process $N$.
In particular, our goal is to find a random set $\hat{\xi}$ that
maximizes the expected value of a specified valuation or gain
function.
The random set $\hat{\xi}$ must be
adapted to the underlying information structure: if the
information available to us at $t\in\mathbf{R}_+^2$ is represented\vspace*{-1pt}
by the
$\sigma$-field $\mathcal{F}_t$, then we must have $\{t\in\hat{\xi}\}\in\mathcal{F}_t$.

\begin{figure}[b]

\includegraphics{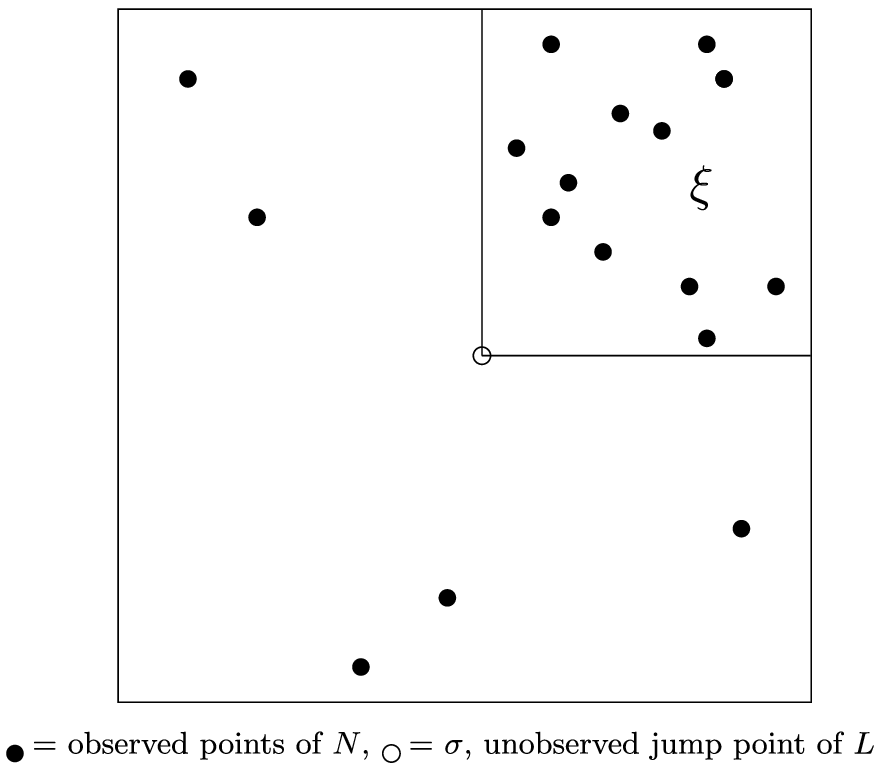}

\caption{A change-set $\xi$ generated by a single point $\sigma$.}\label{Fig1}
\end{figure}

There are many potential areas of application. For example:
\begin{itemize}
\item Environment: The
increased occurrence of polluted wells in a rural area could
indicate a geographic region that has been subjected to industrial
waste.
\item Population health: Unusually frequent outbreaks of a disease
such as leukemia near a nuclear power plant could signal a region of
possible air or ground contamination.
\item Astronomy: A cluster of
black holes could be the result of
an unobservable phenomenon affecting a region in space.
\item Quality control: An increased rate of breakdowns in a
certain type of equipment might follow the failure of one or more
components.
\item Archaeology: An increased number of
archaeological items such as ancient coins found in a particular
region could indicate the location of an event of historical
interest.
\item Forestry: The spread of an airborne disease
through a forest would occur at a higher rate on $\xi$, the set of
points to the northeast of the (unobserved) point ($\sigma$) of
initial infection if the prevailing winds are from the southwest.
\end{itemize}
It is this final type of
example, illustrated in Figure~\ref{Fig1}, that motivates the model to
be studied in this paper.

%

As will be discussed in the conclusion, this paper
represents only a first step in the solution of what we call
the ``optimal set-detection problem.'' Here, we consider the case in
which the change-set $\xi$ is a random \textit{upper layer} (cf.
Section~\ref{Model}) generated by a locally finite set of incomparable
points. In general, the optimal solution $\hat{\xi}$ will be a random
upper layer which is adapted to the available information structure.
This means that the solution is exact in the sense that it is
explicitly defined by the observed data points. This problem cannot be solved
by one-parameter methods. Indeed, even if the random set is
characterized by a single change-point, it will be seen that
the optimal solution does not necessarily correspond to a
point.

In the
one-parameter case, the optimal detection of an exponential change
time in a Poisson process was thoroughly studied in \cite{HJ} using
martingale techniques combined with Bayesian arguments (see also
\cite{PS} for a different approach to the same problem). In
the general set-indexed framework, we found only a very few
papers addressing the problem of a change-point or a
change-set (cf.
\cite{HM}, \cite{KMT1} and \cite{KMT2}). However, none of these
papers deal with the question of the existence of an optimal
solution to the detection problem. Our approach, inspired by that of
\cite{HJ}, makes use of the general
theory of set-indexed martingales as developed in \cite{IM1}. We are
then able to solve the problem with a Bayes-type formula.

The paper is structured as follows. In the next section,
the model is presented and the optimal detection problem is formally
defined. In Section~\ref{Preliminaries}, we give the necessary
background for the multiparameter martingale approach that is the key
for proving the existence of an optimal solution, and develop a
semimartingale representation of the gain function. In Section
\ref{Solution}, sufficient conditions for the existence of
an optimal solution are developed, and then applied to two examples
in Section~\ref{Examples}. Finally, in Section
\ref{Conclusion}, we discuss possible extensions and
directions for further research.

\section{The model}\label{Model}

In order to better understand the two-dimensional model, we
review the change-point problem on $\mathbf{R}_+$ considered in
\cite{HJ}.
We have a nonexplosive point process $N=\{N_t,t\in
\mathbf{R}_+\}$ on $\mathbf{R}_+$, and a random time
$\sigma\geq0$. Given $\sigma$, $N$~is a
Poisson process with intensity $\mu_0$ on $[0,\sigma)$ and
intensity $\mu_1$ on $\xi=[\sigma,\infty)$ ($\mu_1>\mu_0>0$).
Modifying the notation of \cite{HJ} slightly, the gain
function at $t$ is defined by
%
\begin{equation}\label{(1)}
Z_t=c_0(t\wedge\sigma)-c_1(t-\sigma)^++k_0+k_1I_{\{t\geq\sigma\}},
\end{equation}
where $c_0 \geq0$, $c_1>0$ and $k_1\geq0$. The parameters can be
interpreted as follows: the gain function is piecewise linear,
increasing at rate $c_0$ before the jump point and decreasing at rate
$c_1$ after. When $k_1>0$, a penalty equivalent to $-k_1$ is incurred
for stopping the process before the change has occurred. The gain is
maximized when $t=\sigma$.

Let $\mathcal{F}=(\mathcal{F}_t,t\in\mathbf{R}_+)$ denote the
filtration which
characterizes the underlying information available (in
\cite{HJ}, the process $N$ is always $\mathcal{F}$-adapted). For
various filtrations, it is shown in \cite{HJ} that $Z_t$ has a
\textit{smooth semimartingale} (SSM) representation with
respect to $\mathcal{F}$:
%
\begin{equation}\label{(2)}
Z_t=Z_0+\int_0^tU_s\,ds+M_t,
\end{equation}
where $M$ is an $\mathcal{F}$-martingale and $U$ is
$\mathcal{F}$-progressive (i.e., observable). If $U$ is \textit{monotone} in
the sense
that $U_t\leq0\Rightarrow U_{t+h}\leq0 \,\forall h>0$,
then it is straightforward to see that (cf. \cite{HJ}, Theorem
1)
$\hat\sigma:=\inf\{t{}\dvtx{}U_t\leq0\}$ is an optimal $\mathcal{F}$-stopping
rule for $Z$ in terms of expected values: we have
%
\begin{equation}\label{(3)}
E[Z_{\hat\sigma}]=\sup\{E[Z_\tau]{}\dvtx{}\tau\mbox{ an $\mathcal{F}$-stopping time}\}.
\end{equation}
To motivate the model on $\mathbf{R}_+^2$, we will rewrite
(\ref{(1)}) in terms of the single jump point process
$L_t=I_{\{\sigma\leq t\}} $ and the random set
$\xi=[\sigma,\infty)={\{t:L_t>0\}}$:
%
\begin{eqnarray}\label{(4)}
Z_t&=&c_0|A_t\cap\xi^c|-c_1|A_t\cap\xi|+k_0+k_1L_t\nonumber
\\[-8pt]\\[-8pt]
&=& k_0+\int_{A_t} \bigl(-c_1+(c_0+c_1)X_u \bigr)\,du+k_1L_t,\nonumber
\end{eqnarray}
where $A_t=[0,t]$, $|\cdot|$ denotes Lebesgue measure and
$X_t=1-I_{\{t\in\xi\}}=I_{\{L_t=0\}}$.

\begin{figure}[b]

\includegraphics{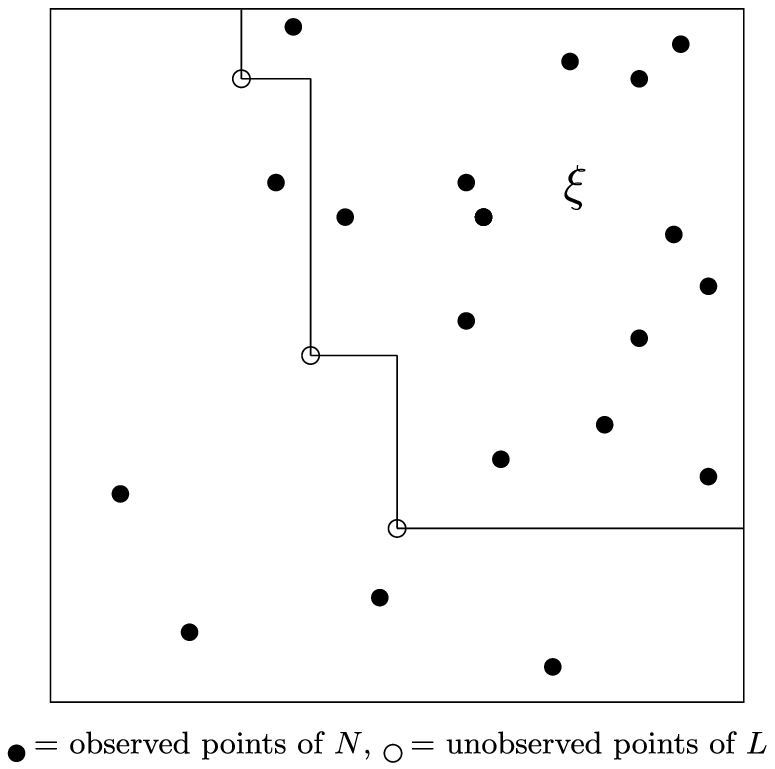}

\caption{A change-set $\xi$ generated by a single line process $L$.}\label{Fig2}
\end{figure}

We are now ready to describe the two-dimensional model. We
are given a random Borel set $\xi\subset(0,\infty)^2$. $N$ is
a nonexplosive point process on $\mathbf{R}^2$ such that given
$\xi$, $N$ is Poisson with intensity $\mu_0$ on $\xi^c$ and
$\mu_1$ on $\xi$. It is always assumed that $\mu_1>\mu_0>0$. (The
case $\mu_0=0$ will be briefly discussed at the end of Section~\ref
{Solution}.) We
will assume that the set $\xi$ is
generated by a
\textit{single line point process}~$L$: that is, $L$ is a
nonexplosive point process whose jump points are all
incomparable ($s,t\in\mathbf{R}_+^2$ are incomparable if
both $s\not\leq t$ and $t\not\leq s$). It is noted
in
\cite{IM2} that in two or more dimensions, the single line process
is the natural generalization of the single jump process, and in
analogy with the change-point model on $\mathbf{R}_+$, we define
$\xi:={\{t{}\dvtx{}L_t>0\}}$. We observe that
$\xi$ is an
\textit{upper layer} ($\xi$ is an upper layer if
$t\in
\xi\Rightarrow s\in
\xi\,\forall s\geq t$). When $L$ has only one
jump point $\sigma$, we observe that
$\xi$ consists of the points to the northeast of $\sigma$. This is
illustrated in Figure~\ref{Fig1}. The more general situation in which
$L$ is a single line process is illustrated in Figure~\ref{Fig2}. In this case, $\xi$ consists of all the points to the northeast of
one or more jump points of~$L$.

%

Using
notation similar to that used for the one-dimensional
problem, for
$t\in \mathbf{R}_+^2$ let
$A_t=\{s\in \mathbf{R}_+^2:s\leq t\}$ and $X_t=1-I_{\{t\in
\xi\}}=I_{\{L_t=0\}}$. The definition of the gain function at
$t\in\mathbf{R}_+^2$ is exactly the same is in (\ref{(4)}):
%
\begin{eqnarray}\label{(5)}
Z_t&=&c_0|A_t\cap\xi^c|-c_1|A_t\cap\xi|+k_0+k_1L_t\nonumber
\\[-8pt]\\[-8pt]
&=&k_0+\int_{A_t} \bigl(-c_1+(c_0+c_1)X_u \bigr)\,du+k_1L_t.\nonumber
\end{eqnarray}
Once again, we assume that $c_0 \geq0$, $c_1>0$ and $k_1\geq
0$, and that $|\cdot|$ denotes Lebesgue measure on $\mathbf{R}_+^2$.

%
\begin{figure}[b]

\includegraphics{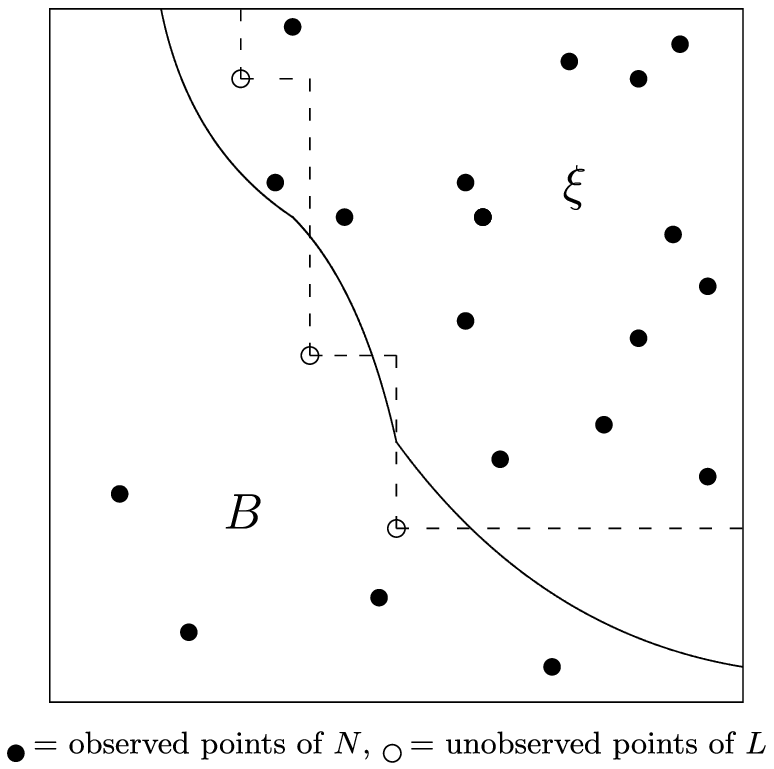}

\caption{A lower layer $B$ and the change-set $\xi$.}\label{Fig3}
\end{figure}

Any point process $N$ can be indexed by the Borel sets in $\mathbf
{R}^2_+$. As in the \hyperref[Introduction]{Introduction}, if $\{\tau_n\}$ denotes the jump
points of $N$ numbered in some arbitrary way, then for any Borel set
$B$, $N(B):=\sum_{n=1}^\infty I(\tau_n\in B)$. [Therefore, we have
$N_t=N(A_t)$.] Consequently, we can define the gain function more
generally over the
class of
\textit{lower layers} $\mathcal{L}$: a set
$B\subseteq\mathbf{R}_+^2$ is a lower layer if $t\in
B\Rightarrow A_t\subseteq B \,\forall t\in\mathbf{R}_+^2$. The gain
function at $B\in\mathcal{L}$ is defined as
%
\begin{eqnarray}\label{(6)}
Z(B)&=&c_0|B\cap\xi^c|-c_1|B \cap\xi
|+k_0+k_1L(B)\nonumber\\[-8pt]\\[-8pt]
&=& k_0+\int_{B} \bigl(-c_1+(c_0+c_1)X_u \bigr)\,du+k_1L(B).\nonumber
\end{eqnarray}
A lower layer $B$ and the change-set $\xi$ are illustrated in Figure~\ref{Fig3};
we observe that $L(B)=2$ in this case.

We see that the gain function defined in (\ref{(6)}) is a natural
generalization of the one-dimensional gain function (\ref{(1)}). The
gain evaluated at $B$ increases in proportion to the area of $B$
outside of the change-set $\xi$, and decreases in proportion to the
area inside of $\xi$. When $k_1>0$, there is a penalty incurred that
is equivalent to $-k_1$ times the number of points in $L$ that lie
outside of (or ``after'') $B$. The gain is maximized when $B=\overline
{\xi^c}$.

We would like to find a random lower layer that maximizes the
expected value of the gain function. The lower layer will
depend on the available information, or more precisely, the
underlying \textit{filtration.}

A class of $\sigma$-fields $\mathcal{F}=\{\mathcal{F}_t,t\in\mathbf
{R}_+^2\}$ is a
filtration if:
\begin{itemize}
\item$\mathcal{F}$ is increasing: $s\leq t\Rightarrow\mathcal{F}_s\subseteq\mathcal{F}_t$, and
\item$\mathcal{F}$ is outer-continuous: $\mathcal{F}_t=\bigcap_{n}\mathcal{F}_{t_n}$ for
every decreasing sequence $(t_n)\subset\mathbf{R}_+^2$ with
$t_n\downarrow t$.
\end{itemize}
%
\begin{Df}[(Cf. \cite{IM1})]\label{dfstopping}
  A closed random lower
layer
$\rho$ is an
$\mathcal{F}$-stopping set if
\[
\{t\in\rho\}\in\mathcal{F}_t \qquad\forall t\in\mathbf{R}_+^2.
\]
\end{Df}

The general optimal set-detection problem in two dimensions
can now be stated as follows: for a given filtration $\mathcal{F}$, our
goal is to maximize
$E[Z_\rho]$, where
$\rho$ is an $\mathcal{F}$-stopping set. If it can be shown that a
stopping set
$\hat\rho$ exists that satisfies the condition
%
\begin{equation}\label{(7)}
E[Z({\hat\rho})]=\sup\{E[Z({\rho})]{}\dvtx{}\rho\mbox{ an $\mathcal{F}$-stopping set}\},
\end{equation}
then our optimal estimate of $\xi$ is
$\hat{\xi}=\overline{\hat\rho^c}$ [$\overline{(\cdot)}$ denotes
set closure]. It is trivial that
$\hat{\xi}$ is an upper layer, and by outer
continuity of $\mathcal{F}$, it is easily seen that $\hat{\xi}$
is also an
adapted random set (i.e.,
$\{t\in\hat{\xi}\}\in\mathcal{F}_t\,\forall t\in\mathbf{R}_+^2)$.

In this paper, we will be focussing on the \textit{sequential}
estimation problem: that is, we will be assuming that
$\mathcal{F}_t=\mathcal{F}^N_t=\sigma\{N_s{}\dvtx{}s\leq t\}$. If $\rho$ is
an $\mathcal{F}
^N$-stopping set, then $I(t\in\rho)$ is a function of the number and
locations of jump points of $N$ in the set $A_t$.
For technical reasons, we shall see that in general it is necessary to
restrict the detection problem to a bounded rectangle
$R=[0,r]^2$. The goal is to find a stopping set
$\hat\rho\subseteq R$ that is optimal in the following sense:
\begin{Df}\label{dfsolution}
 An $\mathcal{F}^N$-stopping set $\hat \rho$
is called an
optimal solution to the sequential detection problem on $R$
provided that $\hat\rho$ satisfies the following equation:
%
\begin{equation}\label{(8)}
E[Z({\hat\rho})]=\sup\{E[Z({\rho})]{}\dvtx{}\rho\subseteq R
\mbox{ an $\mathcal{F}^N$-stopping set}\}.
\end{equation}
\end{Df}

Restricting our attention to $R$ ensures that $\hat\rho$
is bounded and so $E[Z({\hat\rho})]$ is always well defined.
In this case, we have an optimal estimate $\hat{\xi}_R $ of
$\xi\cap R$, defined by
$\hat{\xi}_R=\overline{R\setminus\hat\rho}$.

\section{Mathematical preliminaries}\label{Preliminaries}

In this section we present the mathematical tools needed in the
sequel. In \cite{HJ}, Herberts and Jensen make
use of martingale techniques to provide a simple and elegant
method of finding sufficient conditions for the existence of an
optimal solution to the detection problem on $\mathbf{R}_+$.
Martingale methods have been extended to more general spaces
in \cite{IM1}, and we are able to exploit this theory in a
similar way. To motivate the necessary technical details that follow,
we first describe our overall plan of attack. Recall that $\mathcal{F}^N$
denotes the filtration representing the data that can be observed, and
below $\mathcal{G}$ will denote a larger filtration containing additional
information, some of which cannot be observed.

\textit{Plan of attack}:
\begin{itemize}
\item The gain function $Z$ can be rewritten as a (two-parameter)
semimartingale (Definition~\ref{defssm}):
\[
Z_B=k_0+\int_BU_t\, dt +k_1M_B,
\]
where $M$ is a weak martingale (Definition~\ref{mg}) with respect to a
filtration $\mathcal{G}$ and $U$ is $\mathcal{G}$-adapted but not necessarily
observable (cf. Lemma~\ref{fullssm}).
%
\item For the observable filtration $\mathcal{F}^N$ and $\rho$ an
$\mathcal{F}
^N$-stopping set, we have $E[M_\rho]=0$ (Lemma~\ref{stopping}) and if
$V_t=E[U_t|\mathcal{F}_t^N]$ (observable), then\break (Lemma~\ref{projection})
\[
E[Z_\rho]=k_0+E \biggl[\int_\rho U_t\,dt \biggr]\\
= k_0+E \biggl[\int_\rho V_t\,dt \biggr].
\]
\item If $V$ satisfies a monotonicity property on $R$ (cf.
Definition~\ref{defmonotone} and\break Lemma~\ref{projection}), then there exists
an optimal solution $\hat\rho$ to
the sequential detection problem on~$R$, and
the optimal estimate of $\xi\cap R$ is
\[
\hat\xi_R=\{t\in R{}\dvtx{} V_t\leq0\}.
\]
\end{itemize}

Keeping this outline of our approach in mind, we continue with the
necessary mathematical details.

\subsection{Martingale preliminaries}\label{Martingale}

Martingales on $\mathbf{R}_+^2$ can be defined in various
ways (cf. \cite{IM1}), but here we need only the weakest
definition. In what follows, $T$~denotes either $\mathbf{R}_+^2$ or a
bounded region $R=[0,r]^2$, and $(\Omega,
\mathcal{F}, P)$ is a complete probability space equipped with a
$T$-indexed filtration $\mathcal{F}=\{\mathcal{F}_t{}\dvtx{}t\in T\}$
(without loss of
generality, assume that $\mathcal{F}_t$ contains all the $P$-null
sets $\forall t\in T$). A~$T$-indexed process
$X=\{X_t{}\dvtx{}t\in T\}$ is adapted to $\mathcal{F}$ if $X_t$ is $\mathcal
{F}_t$-measurable,
for all $ t\in T$. For any
$T$-indexed process
$X=\{X_t{}\dvtx{}t\in T\}$, for $s=(s_1,s_2)\leq(t_1,t_2)=t\in T$,
define the increment of $X$ over the rectangle
$(s,t]=(s_1,t_1]\times(s_2,t_2]$ in the usual way:
\[
X(s,t]=X_{(t_1,t_2)}-X_{(s_1,t_2)}-X_{(t_1,s_2)}+X_{(s_1,s_2)}.
\]
%
\begin{Df}\label{mg}
 Let $M=\{M_t{}\dvtx{}t\in T\}$ be an integrable
process on
$T$, adapted to a filtration $\mathcal{F}=\{\mathcal{F}_t{}\dvtx{}t\in T\}$.
$M$ is a weak
$\mathcal{F}$-(sub)martingale if $M$ is equal to 0 on the axes, and for
every
$s\leq t\in T$,
\[
E[M(s,t]|\mathcal{F}_s]=(\geq)0.
\]
(A process
$X$ is integrable if $E[|X_t|]<\infty \,\forall t\in
T$.)
\end{Df}
%
\begin{Df}\label{increasing}
Let $v=\{v_t{}\dvtx{}t\in T\}$ be a
function on
$T$. We say that $v$ is increasing (decreasing) if:
\begin{itemize}
\item$v$ is 0 on the axes,
\item$v$ is outer continuous with inner limits:
that is, $v$ is continuous from above and with limits from the
other three quadrants at each
$t\in T$, and
\item for every $s\leq t\in T$, $v(s,t]\geq(\leq) 0$.
\end{itemize}
A process $V=\{V_t{}\dvtx{}t\in T\}$ is increasing (decreasing) if
for each $\omega\in\Omega$, the function
$ V_{\bolds{\cdot}}(\omega)$ is increasing (decreasing).
\end{Df}

\begin{Com}\label{measure}
An increasing function $v$ can be regarded as the distribution
of a measure on $\mathbf{R}_+^2$. Therefore, $v(B)$ is well defined
for any
Borel set
$B$, where we use $v_{\bolds{\cdot}}$
and $v(\cdot)$ to denote, respectively, the function and the
generated measure. Likewise, a decreasing function generates
a negative measure, and we will use similar notation.
\end{Com}
%
\begin{Df}{ Let $L$ be a weak
$\mathcal{F}$-submartingale. An increasing process $\Lambda$ is a
compensator for
$L$ if
$\Lambda$ is
$\mathcal{F}$-adapted and $M=L-\Lambda$ is a weak martingale. }
\end{Df}

\begin{Com}
 As defined above, the compensator of a submartingale
need not be unique (any increasing process is trivially a compensator
for itself). A type of predictability is required for uniqueness (cf.
\cite{IM1}), but this point is not of importance here.
\end{Com}

In light of Comment~\ref{measure}, the following lemma is a
special case of of Lemma~3.3.5 of
\cite{IM1}.
%
\begin{Lem}\label{stopping}
 If $M$ is a weak martingale which can be
expressed as the difference of two increasing integrable processes,
and
$\rho$ is a stopping set such that $\rho\subseteq R=[0,r]^2$,
then
$M({\rho})$ is well defined and $E[M({\rho})]=0$.
\end{Lem}
%
\begin{Df}\label{defssm}
Let $Z=\{Z_t{}\dvtx{}t\in T\}$ be a process on $T$,
adapted to a filtration $\mathcal{F}=\{\mathcal{F}_t{}\dvtx{}t\in T\}$. $Z$
is a
smooth semimartingale with respect to $\mathcal{F}$ ($\mathcal{F}$-SSM)
if it
satisfies a decomposition of the form
%
\begin{equation}\label{(9)}
Z_t=Z_{(0,0)}+\int_0^{t_1}\int_0^{t_2}U_{(s_1,s_2)}\,ds_2\,ds_1+M_t
\end{equation}
for each $t=(t_1,t_2)\in T$, where $U$ is an outer continuous
process with inner limits adapted to $\mathcal{F}$ and $M$ is a weak
$\mathcal{F}$-martingale. We denote the $\mathcal{F}$-SSM as $Z=(U,M)$.
\end{Df}

In order to show that an optimal solution exists to the sequential
detection problem, we will require a monotonicity
property.
\begin{Df}\label{defmonotone}
A function $v=\{v_t{}\dvtx{}t\in T\}$ is monotone on $T$ if $v_s\leq
0\Rightarrow v_t\leq0\,\forall t\geq s\in T$. A
process
$V$ is monotone if $V_{\bolds{\cdot}}(\omega)$ is monotone for each
$\omega\in\Omega$.
\end{Df}

\begin{Com}
\begin{enumerate}
\item Note that any decreasing function is monotone, but the
converse is not true.
\item If a process $V$ is decreasing in each parameter
separately on $T$, then $V$ is monotone on $T$ but not necessarily
decreasing in the sense of Definition~\ref{increasing}.
\item Note that if $V$ is monotone, then $V_t>0\Rightarrow
V_s>0\,\forall s\leq t$.
\item If $V$ is monotone and adapted to a filtration $\mathcal{F}$, the
set
%
\begin{equation}\label{(10)}
\hat\rho=\{t\in T{}\dvtx{} V_s>0\,\forall s\ll t\}
\end{equation}
is an $\mathcal{F}$-stopping set (cf. Definition~\ref{dfstopping}).
[$s=(s_1,s_2)\ll(t_1,t_2)=t
\Leftrightarrow s_i<t_i$ if $t_i>0$, and
$s_i=0$ if $t_i=0$, $i=1,2$.] Clearly, $\hat\rho$ is a random closed
lower layer, and the fact that $V$ is adapted ensures that
$\{t\in\hat\rho\}\in\mathcal{F}_t$: taking any sequence
$(t_n) \uparrow t$ with $t_n\ll t$, by monotonicity it follows
that
\[
\{t\in\hat\rho\}=\bigcap_n\{V_{t_n}>0\}\in
\bigcup_n\mathcal{F}_{t_n}\subseteq\mathcal{F}_t.
\]
\end{enumerate}
\end{Com}

In \cite{HJ}, the solution to the optimal stopping problem is based on
a SSM representation of the form (\ref{(2)}), which in turn is based
on a projection theorem. The question of the existence of optional
and predictable projections in higher dimensions is a delicate one,
usually requiring a strong assumption of conditional independence
on the underlying
filtration [denoted (F4) in the two-dimensional literature]. For
details, see \cite{MZ}, for example. In practice, one can generally
show directly that a suitable projection exists without relying on a
general existence theorem, and for our purposes the following lemma
will be adequate.
\begin{Lem}\label{projection}
Let $U$ be a bounded $T$-indexed process adapted to a
filtration
$\mathcal{G}$ such that
$U$ is outer-continuous
with inner limits. If
$\mathcal{F}$ is a subfiltration of $\mathcal{G}$ (i.e., $\mathcal{F}_t\subseteq\mathcal{G}_t\,\forall
t$), and if a version of $ V_t=E[U_t|\mathcal{F}_t]$ exists that is
outer-continuous with inner limits, then for any $\mathcal
{F}$-stopping set
$\rho\subseteq R=[0,r]^2$,
%
\begin{equation}\label{(10a)}
E \biggl[\int_\rho U_t\,dt \biggr]=E \biggl[\int_\rho V_t\,dt \biggr].
\end{equation}
In addition, if $V$ is monotone on $R$, then the $\mathcal
{F}$-stopping set
$\hat\rho\subseteq R$ defined by
%
\begin{equation}\label{(11)}
\hat\rho=\{t\in R{}\dvtx{} V_s>0\,\forall s\ll t\}
\end{equation}
is optimal in the sense that
\[
E \biggl[\int_{\hat\rho}U_t\,dt \biggr]
=\sup \biggl\{ E \biggl[\int_\rho U_t\,dt \biggr]{}\dvtx{}\rho\subseteq R, \rho\mbox{ an $\mathcal{F}$-stopping
set}  \biggr\}.
\]
\end{Lem}

\begin{pf}
First, the assumption that $U$ and $V$ have sample
paths that are regular (outer-continuous with inner limits) and that
$U$ (and hence $V$) is bounded ensures that the integrals and
expectations in (\ref{(10a)}) are well defined.

Next, let $T_n:=\{(\frac{i }{2^n}r,\frac{j }{2^n}r){}\dvtx{} 0\leq i,j\leq2^n\}$ denote the ``dyadics'' of order $n$ in~$R$. The class of
rectangles $\mathcal{C}_n $ partitions $R$, where $C\in\mathcal{C}_n$ if
$C$ is of the form $C=A_t\setminus(\bigcup_{s\in T_n, s\not\geq t}A_s)$
for some $t\in T_n$. Let $t_{C-}=\inf\{t\in C\}$\vspace*{1pt} denote the lower
left corner of $C$. We now define the ``discrete'' approximation
$\rho_n$ of $\rho$ by\vspace*{-1pt}
\[
\rho_n=\bigcup_{C\in\mathcal{C}_n:t_{C-}\in\rho}C.
\]
It is straightforward that $\rho_n\subseteq R$ is an $\mathcal{F}$-stopping
set, that
$(\rho_n)$ is decreasing in $n$ and $\rho=\bigcap_n \rho_n$.
Boundedness and uniform integrability ensure that
$E [\int_\rho U_t\,dt ]=\lim_nE [\int_{\rho_n} U_t\,dt ]$ and
$E [\int_\rho V_t\,dt ]=\lim_nE [\int_{\rho_n} V_t\,dt ]$.
To complete\vspace*{1pt} the proof of the first statement in the theorem, observe
that by boundedness of $U$,
\begin{eqnarray*}
E \biggl[\int_{\rho_n}U_t\,dt \biggr]
&=&E \biggl[\sum_{C\in\mathcal{C}_n}I_{\{t_{C-}\in\rho\}}\int_C U_t \,dt \biggr]\\
&=&E \biggl[\sum_{C\in\mathcal{C}_n}I_{\{t_{C-}\in\rho\}}E \biggl[\int_C U_t \,dt\Big|\mathcal{F}_{t_{C-}} \biggr] \biggr]\\
&=&E \biggl[\sum_{C\in\mathcal{C}_n}I_{\{t_{C-}\in\rho\}}E \biggl[\int_C E[U_t|\mathcal{F}_t] \,dt\Big|\mathcal{F}_{t_{C-}} \biggr] \biggr]\\
&=&E \biggl[\sum_{C\in\mathcal{C}_n}I_{\{t_{C-}\in\rho\}}E \biggl[\int_C V_t \,dt\Big|\mathcal{F}_{t_{C-}} \biggr] \biggr]\\
&=&E \biggl[\sum_{C\in\mathcal{C}_n}I_{\{t_{C-}\in\rho\}}\int_C V_t \,dt \biggr]
=E \biggl[\int_{\rho_n}V_t\,dt \biggr].
\end{eqnarray*}
The third equality above follows by Fubini and the assumption
that $V$ has regular sample paths, and since
$t\in C\Rightarrow t\geq t_{C-}$. [The assumption that $V$ has a
version with regular sample
paths ensures that $V$ is jointly $\mathcal{F}\times\mathcal
{B}(\mathbf
{R}_+^2)$-measurable,
where
$\mathcal{B}(\mathbf{R}_+^2)$ denotes the Borel sets in $R$.]

Next, assume that $V$ is monotone. To prove optimality of $\hat\rho$,
let
$\rho\subseteq R$ be any other stopping set in $R$. We have
\[
E \biggl[ \int_{\hat\rho} U_t \,dt -\int_\rho U_t\,dt \biggr]
=E \biggl[ \int_{\hat\rho\setminus\rho} V_t \,dt
-\int_{\rho\setminus\hat\rho} V_t\,dt \biggr]\geq0,
\]
since $V>0$ on $\hat\rho^o$ (the interior of $\hat\rho$) and $V\leq0$
on
$\hat\rho^c$.
\end{pf}

\subsection{Smooth semimartingale representation of the gain function}\label{SSM}

We begin this section with an analysis of the single line process
$L_t$:
$L$ is a nonexplosive point process whose jump points are all
incomparable. Single line processes and their compensators were
discussed
in
\cite{IM2}, to which the reader may refer for more detail.
Heuristically, if $\mathcal{F}^L_s=\sigma(L_u{}\dvtx{}u\leq s)$, then a process
$\Lambda$ will be an
$\mathcal{F}^L$-compensator
of
$L$ if
\begin{eqnarray*}
&&\Lambda((s_1,s_2),(s_1+ds_1,s_2+ds_2)])
\\
&&\qquad\approx I_{\{L_s=0\}}E\bigl[ L ((s_1,s_2),(s_1+ds_1,s_2+ds_2)]) |L_s=0\bigr],
\end{eqnarray*}
since $L$ cannot have any jump points in
$((s_1,s_2),(s_1+ds_1,s_2+ds_2)])$ if $L_s>0$ and $\{L_s=0\}$
is an atom of $\mathcal{F}^L_s$. Define the (deterministic) increasing
function
$\Lambda^{(s)}_t:=E[L( s,t]|L_s=0]$, for $t\geq s$, and when it
exists, let
\[
\lambda_s=\lim_{t_1\downarrow s_1,t_2\downarrow
s_2}\frac{\Lambda^{(s)}_t }{(t_1-s_1)(t_2-s_2)}.
\]
In particular, if
$\lambda_s$ exists
for every $s\in T$ and is Lebesgue measurable,
then
%
\begin{eqnarray}\label{(12)}
\Lambda_t=\int_{A_t}\lambda_uI_{\{L_u=0\}}\,du.
\end{eqnarray}

In what follows (and as will be seen to be the case in our examples),
we will assume that a representation of the form (\ref{(12)}) exists
for the compensator $\Lambda$ of $L$, and we will refer to the
deterministic function $\lambda$ as the \textit{weak hazard function}
of $L$. It will always be assumed that $\lambda$ is continuous.

To better understand the weak hazard, we observe
that if\vspace*{1pt}
$E[L]$ of $L$ is absolutely continuous with respect to Lebesgue measure
with Radon--Nikodym
derivative $\tilde\lambda$, then for every $u\in T$ with
$P(L_u=0)>0$,
$\lambda_u=\tilde\lambda_u/P(L_u=0)$. To see this, simply
observe that for each $ t\in T$,
%
\begin{equation}\label{(13)}
\int_{A_t}\tilde\lambda_u\,du=E[L_t]=E[\Lambda_t]=\int_{A_t}
\lambda_uP(L_u=0)\,du.
\end{equation}

Returning to the gain function (\ref{(6)}), let $M$ denote
the weak martingale $L-\Lambda$ and recall that $X_u=I_{\{L_u=0\}}$.
For any lower layer $B\subseteq T$,
%
\begin{eqnarray}\label{(14)}
Z(B)&=& k_0+\int_{B} \bigl(-c_1+(c_0+c_1)X_u \bigr)\,du+k_1L(B)\nonumber
\\[-8pt]\\[-8pt]
&=& k_0+\int_{B} \bigl(-c_1+(c_0+c_1+k_1\lambda_u)X_u \bigr)\,du+k_1M(B).\nonumber
\end{eqnarray}
We note that $X$ is outer-continuous with inner limits by
definition and that $\lambda$ is assumed to be continuous, and so we
now have an
$\mathcal{F}^L$-SSM representation of the gain function:
$Z=(U,M)$, where $U_t:= -c_1+(c_0+c_1+k_1\lambda_t)X_t $.

\begin{Com}\label{noinfo}
As a simple illustration, if the point
process
$L$ and the set
$\xi=\{t{}\dvtx{}L_t>0\}$ are unobservable and no other information is
available (i.e., $N$ is not observed and
$\mathcal{F}_t=\{\varnothing,\Omega\} \,\forall t\in T$), then for
$R=[0,r]^2$, we
are looking for a deterministic set
$\hat B\subseteq R$ that maximizes\vspace*{-2pt}
%
\begin{eqnarray}
E[Z(B)]&=&E \biggl[k_0+\int_{B} \bigl(-c_1+(c_0+c_1+k_1\lambda_u)X_u \bigr)\, du+k_1M(B)\biggr]\nonumber\\[-8.5pt]\\[-8.5pt]
&=& k_0+\int_{B} \bigl(-c_1+(c_0+c_1+k_1\lambda_u )P(L_u=0)\bigr)\,du.\nonumber
\end{eqnarray}
Letting $V_u=[-c_1+(c_0+c_1+k_1\lambda_u )P(L_u=0)]$, it is easily
seen that $V$ is deterministic and an
optimal solution for the detection problem exists if $V$ is monotone,
in which case\vspace*{-2pt}
%
\begin{eqnarray}\label{(16)}
\hat B &=&\{t\in R{}\dvtx{}V_u> 0 \,\forall u\ll t\}\nonumber\\[-9pt]\\[-9pt]
&=& \biggl\{t\in R{}\dvtx{}P(L_u=0)>\frac{c_1}{(c_0+c_1+k_1\lambda_u )}\,\forall u\ll t \biggr\}.\nonumber
\end{eqnarray}
The optimal estimate of $\xi\cap R$ is
\[
\hat{\xi}_R= \biggl\{t\in R{}\dvtx{}P(L_t=0)\leq\frac{c_1}{(c_0+c_1+k_1\lambda_t )} \biggr\}.
\]
\end{Com}

\begin{Ex}[(The single jump process)]\label{singlejump}
Suppose $L_t=I_{\{Y\in A_t\}}$, where $Y$ is a $T$-valued random
variable with distribution $F$ and continuous density $f$. Then we
have
$\lambda_u=\frac{f_u}{1-F_u}$. To verify that the representation
(\ref{(12)})\vspace*{-1pt} is satisfied with this definition, observe first that
$E[L(s,t]|\mathcal{F}_s^L]=\frac{F(s,t]}{1-F_s}I_{\{L_s=0\}}.$ Next,\vspace*{-2pt}
\begin{eqnarray*}
E \biggl[\int_{(s,t]}\frac{f_u}{1-F_u}I_{\{L_u=0\}}\,du\Big|\mathcal{F}_s \biggr]
&=&\int_{(s,t]}\frac{f_u}{1-F_u}P(L_u=0|\mathcal{F}_s)\,du\\[-1pt]
&=&\int_{(s,t]}\frac{f_u}{1-F_u}\cdot\frac{1-F_u}{1-F_s}I_{\{L_s=0\}}\,du\\[-1pt]
&=&\frac{F(s,t]}{1-F_s}I_{\{L_s=0\}}.
\end{eqnarray*}
Thus, the increasing process
$\Lambda_t=\int_{A_t}\lambda_sI(L_s=0)\,ds$ is a $\mathcal
{F}^L$-compensator for
$L$, verifying (\ref{(12)}).

It should be noted that in the literature on bivariate survival
analysis, the definition of the hazard function is
$\frac{f_u}{S_u}$ where $S_u=P(Y\geq u)$. For this reason, we refer
to our hazard $\lambda=\frac{f}{1-F}$ as the ``weak'' hazard.

Returning to Comment~\ref{noinfo}, when no information is available,
$V$ is decreasing and
(\ref{(16)}) defines an optimal deterministic solution if $f$
is decreasing in each parameter.
\end{Ex}

\begin{Ex}[(First line of a Poisson process)]\label{poisson}
Consider a homogeneous Poisson process $J$ on $T$ with rate $\gamma$.
If
$\Delta_J$ denotes the set of jump points of $J$, then the
\textit{first line} of $J$ is the single line point process $L$ with
(incomparable) jump points
\[
\Delta_L=\min(\Delta_J)=\{\tau\in
\Delta_J{}\dvtx{}\tau'\not\leq\tau\,\forall\tau'\in\Delta_J\mbox{ such
that }\tau'\neq\tau\}.
\]
In this case,
$\xi=\{t{}\dvtx{}L_t>0\}=\{t{}\dvtx{}J_t>0\}$. As is shown in \cite{IM2}, the
weak hazard of $L$ is $\gamma$.

Considering the situation in Comment~\ref{noinfo} when no
information is available, we have $V_u=-c_1+(c_0+c_1+k_1\gamma
)e^{-\gamma
u_1u_2} $, which is clearly monotone. In this case, the optimal
solution given in (\ref{(16)}) becomes
\begin{eqnarray*}
\hat B&=& \biggl\{t\in R{}\dvtx{}e^{-\gamma t_1t_2}\geq\frac{c_1}{(c_0+c_1+k_1\gamma)} \biggr\} \\
&=& \biggl\{t\in R{}\dvtx{} t_1t_2 \leq\frac{\ln(c_0+c_1+k_1\gamma)-\ln(c_1)}{\gamma} \biggr\}
\end{eqnarray*}
and
\[
\hat\xi_R= \biggl\{t\in R{}\dvtx{} t_1t_2 \geq\frac{\ln(c_0+c_1+k_1\gamma)-\ln(c_1)}{\gamma} \biggr\}.
\]
\end{Ex}

We are now ready to return to the sequential detection problem, and
consider the case in which the process
$N$ is observed (recall that $N$~is a Poisson process with rate
$\mu_0$ on
$\xi^c$ and
$\mu_1$ on $\xi$). We denote the full filtration
$\mathcal{F}^{L,N}=\{\mathcal{F}^{L,N}_t{}\dvtx{}t\in T\}$, where
$\mathcal{F}^{L,N}_t=\sigma\{L_s,N_s, s\leq t\}$, and (as before) the
subfiltrations
$\mathcal{F}^L=\{\mathcal{F}^L_t{}\dvtx{}t\in T\} $ and $\mathcal{F}^N=\{
\mathcal{F}^N_t{}\dvtx{}t\in T\}$ where
$\mathcal{F}^L_t=\sigma\{L_s{}\dvtx{}s\leq t\}$ and $\mathcal{F}^N_t=\sigma
\{N_s{}\dvtx{}s\leq t\}$.
Although we defined the weak hazard of
$L$ with respect to $\mathcal{F}^L$, it is easy to see that given the
full filtration
$\mathcal{F}^{L,N}$,
$L-\Lambda$ is still a weak $\mathcal{F}^{L,N}$-martingale. This
follows because on
$\{L_s=0\}=\{s\in\xi^c\}$,
$N$ is a Poisson process with rate $\mu_0$ on $A_s$ and so $N|_{A_s}$
($N$ restricted to $A_s$) adds no additional information about the
behavior of
$L_t$ for
$t>s$. Formally, we have
\[
E[L(s,t]|\mathcal{F}^{L,N}_s]=I_{\{L_s=0\}}\Lambda
^{(s)}_t=E[L(s,t]|\mathcal{F}^{L}_s].
\]
Therefore, from this discussion we have the following lemma
and we are ready to proceed with finding an optimal solution
to the sequential detection problem.
%
\begin{Lem}\label{fullssm}
Equation (\ref{(14)}) defines an $\mathcal{F}^{L,N}$-SSM representation
of the gain
function $Z{}\dvt{} Z=(U,M)$ where $U_t:= -c_1+(c_0+c_1+k_1\lambda_t)X_t $.
\end{Lem}

\section{Optimal solution to the sequential detection problem}\label{Solution}

We consider the $\mathcal{F}^{L,N}$-SSM representation of the gain function
(\ref{(14)}):
\[
Z(B)= k_0+\int_{B} \bigl(-c_1+(c_0+c_1+k_1\lambda_u)X_u \bigr)\,du+k_1M(B).
\]
In order to find sufficient conditions for the existence of an optimal
solution in the sequential case, we will be appealing to Lemma
\ref{projection}, with
$\mathcal{G}=\mathcal{F}^{L,N}$,
$\mathcal{F}=\mathcal{F}^N$ and $U_t= -c_1+(c_0+c_1+k_1\lambda
_t)X_t$. In
order to find $V_t=E[U_t|\mathcal{F}^N_t]$, it is enough to determine
\[
E[X_t|\mathcal{F}^N_t]=P(L_t=0|\mathcal{F}_t^N).
\]
As in \cite{HJ},
we use a Bayesian argument. The first step is to determine the
conditional likelihood $\ell_{N|L}(t)$
of $N|_{A_t}$ 
given
$L$ and use this to find the likelihood $\ell_N(t)$ of
$N|_{A_t}$. Next we find the conditional likelihood $\ell_{N|L_t=0}(t)$
of
$N|_{A_t}$ on the set
$\{L_t=0\}$. Finally, we have
%
\begin{eqnarray}\label{(17)}
E[X_t|\mathcal{F}^N_t]&=&P(L_t=0|\mathcal{F}_t^N)\nonumber\\[-8pt]\\[-8pt]
&=&\frac{\ell_{N|L_t=0}(t)\times P(L_t=0)}{\ell_N(t)}.\nonumber
\end{eqnarray}

When computing the likelihood $\ell_{N|L}$, in fact it is equivalent to
condition on the random {upper layer} $\xi=\{u{}\dvtx{}L_u>0\}$. To see
this, let
$(\mathcal{U},d_H)$ denote the collection of closed upper layers in
$T$ endowed
with the Hausdorff metric. It is shown in
\cite{IMP} that $(\mathcal{U}, d_H)$ is a complete separable metric space
and that $\xi$ can be regarded as the unique jump point in a single
jump process $\tilde L$ on $\mathcal{U}$; in addition, $L$ determines
and is
determined by $\tilde L$. In particular, $L_t>0\Leftrightarrow t
\in
\xi\Leftrightarrow E_t \subseteq\xi$, where
$E_t=\{s\in T{}\dvtx{}s\geq t\}$. Let $\mu_\xi$ denote the
measure induced by $\xi$ on
$\mathcal{U}$.

Given $L$, or equivalently $\xi$, $N$ is a Poisson process with rate
$\mu_0$ on $\xi^c$ and $\mu_1$ on $\xi$. Using the well-known
likelihood for the Poisson process (cf. \cite{DV}, page~22), we have
%
\begin{eqnarray}\label{(18)}
\ell_{N|L}(t)&=&\ell_{N|\xi}(t)\nonumber\\
&=&e^{-\mu_0|A_t\setminus\xi|}\mu_0^{N(A_t\setminus\xi)}e^{-\mu_1|A_t\cap\xi|}\mu_1^{N(A_t\cap\xi)}\\
&=&e^{-\mu_0|A_t|}\mu_0^{N_t}e^{-(\mu_1-\mu_0)|A_t\cap\xi|} \biggl(\frac{\mu_1}{\mu_0} \biggr)^{N(A_t\cap\xi)}.\nonumber
\end{eqnarray}
By considering separately the events $\{L_t=0\}=\{t\notin
\xi\}=\{E_t\not\subseteq\xi\}=\{A_t\cap\xi=\varnothing\}$ and $\{
L_t>0\}=\{t\in
\xi\}=\{E_t\subseteq
\xi\} $, we use
(\ref{(18)}) obtain
%
\begin{eqnarray}\label{(19)}
\ell_N(t)&=&P(L_t=0)e^{-\mu_0|A_t|}\mu_0^{N_t}\nonumber
\\ &&+e^{-\mu_0|A_t|}\mu_0^{N_t} {\int_{\{D\in\mathcal{U}:E_t\subseteq D\}}}e^ {-(\mu_1-\mu_0)|A_t\cap D|}
\biggl(\frac{\mu_1}{\mu_0} \biggr)^{N(A_t\cap D)}\,d\mu_\xi(D)\hspace*{-28pt}\\
&=&e^{-\mu_0|A_t|}\mu_0^{N_t} \bigl[P(L_t=0)+e^{-(\mu_1-\mu_0)|A_t|}Q_t \bigr],\nonumber
\end{eqnarray}
%
where
\begin{equation}\label{(19a)}
Q_t= {\int_{\{D\in\mathcal{U}:E_t\subseteq D\}}}e^ {(\mu_1-\mu_0)|A_t\setminus
D|} \biggl(\frac{\mu_1} {\mu_0} \biggr)^{N(A_t\cap D)}\,d\mu_\xi(D).
\end{equation}

Before continuing, we observe that since $\mu_1>\mu_0$, $Q$ is
increasing in each parameter separately because
each term in the integrand is increasing in each component
for
$D$ fixed, and the range of integration is increasing since the set
$E_t$ decreases with each component.

Next, if $L_t=0$, then $N|_{A_t}$ is Poisson with rate $\mu_0$, and
so
%
\begin{equation}\label{(20)}
\ell_{N|L_t=0}=e^{-\mu_0|A_t|}\mu_0^{N_t}.
\end{equation}
Substituting (\ref{(19)})
and (\ref{(20)}) in (\ref{(17)}), we obtain
%
\begin{eqnarray}\label{(21)}
E[X_t| \mathcal{F}_t^N]&=&\frac{e^{-\mu_0|A_t|}\mu_0^{N_t}P(L_t=0)}
{e^{-\mu_0|A_t|}\mu_0^{N_t} [P(L_t=0)+e^{-(\mu_1-\mu_0)|A_t|}
Q_t ]}\nonumber\\[-8pt]\\[-8pt]
&=&\frac{1}{1+q_tQ_t},\nonumber
\end{eqnarray}
where $q_t=\frac{e^{-(\mu_1-\mu_0)|A_t|}}{P(L_t=0)}$. [If $
P(L_t=0)=0$, (\ref{(21)}) remains formally valid since $E[X_t|
\mathcal{
F}_t^N]=0
$ and
$q_t=\infty$.] We are now ready to state our main
result:
\begin{Th}\label{main}
 Let $L$ be a single line process with
continuous weak hazard~$\lambda$, and define the function $q$ by
\[
q_t=\frac{e^{-(\mu_1-\mu_0)t_1t_2}}{P(L_t=0)}\qquad \mbox{for }
t=(t_1,t_2)\in\mathbf{R}_+^2.
\]
An optimal solution to the sequential
detection problem on
$R=[0,r]^2$ exists if
$\lambda$ and
$q$ are decreasing and increasing, respectively, in each component on
$R$. In this case $V$ is monotone on
$R$, and the optimal solution is given by (\ref{(11)}):
\[
\hat\rho=\{t\in R{}\dvtx{} V_s>0\,\forall s\ll t\},
\]
where
\[
V_t=-c_1+(c_0+c_1+k_1\lambda_t)\frac{1}{1+q_tQ_t}.
\]
\end{Th}
\begin{pf}
We review our results so far. We have
the $\mathcal{F}^{L,N}$-SSM representation of the gain function
$Z(B)= k_0+\int_BU_t\,dt+k_1M(B)$, where
$U_t= -c_1+(c_0+c_1+k_1\lambda_t)X_t $. $U$ is bounded since
$\lambda$ is decreasing in each component and $X$ is an indicator
function. By the argument immediately preceding the theorem, we have
that
%
\begin{equation}\label{(22)}
V_t=E[U_t|\mathcal{F}^N_t]=-c_1+(c_0+c_1+k_1\lambda_t)\frac{1}{1+q_tQ_t}.
\end{equation}
To see that $V$ has a version which is outer-continuous with inner
limits (o.c.i.l.), recall that $\lambda$ is assumed to be continuous and
observe that
$q$ is o.c.i.l. by definition. Turning next to $Q$, we see that the
integrand in (\ref{(19a)}) is o.c.i.l. and increasing in each
component in
$t$, as is
\[
\mu_\xi(\{D\in\mathcal{U}{}\dvtx{}E_t\subseteq D\})=P(L_t>0).
\]
Therefore, it follows that $Q$, and hence $V$ are o.c.i.l. Therefore,
Lemmas~\ref{stopping} and~\ref{projection} imply that for any
$\mathcal{F}^N$-stopping set
$\rho\subseteq R$,
%
\begin{eqnarray}
E[Z({\rho})]
= k_0+E \biggl[\int_\rho U_t\,dt \biggr]
= k_0+E \biggl[\int_\rho V_t\,dt \biggr].
\end{eqnarray}
To show that an optimal solution $\hat\rho$ exists [as in
(\ref{(11)})], it is sufficient to show that $V$ is monotone (again,
by Lemma~\ref{projection}). Since we have already seen that $Q$ is
increasing in each component on $R$, the assumption that
$\lambda$ and $q$ are decreasing and increasing, respectively, in
each component imply that
$V$ is monotone on $R$.
This completes the proof.
\end{pf}

\begin{Com}
It has been pointed out by an anonymous referee that
the case $\mu_0=0$ relates to a so-called support estimation problem.
In this case, the random set $\xi$ denotes the support of a Poisson
process with rate $\mu_1$. The gain function can be defined exactly as
before, and the analysis proceeds in very much the same way. Now we
know that $N_t>0\Rightarrow t\in\xi\Rightarrow L_t>0$, and equation~(\ref{(17)}) becomes
%
\begin{eqnarray}\label{(17a)}
E[X_t|\mathcal{F}^N_t]
&=&P(L_t=0|\mathcal{F}_t^N)\nonumber\\
&=& P(L_t=0|N_t=0)I(N_t=0)\\
&=&\frac{P(L_t=0)}{P(N_t=0)}I(N_t=0).\nonumber
\end{eqnarray}
%
Continuing with the same sort of arguments used previously, if $\mu
_0=0$, equation~(\ref{(21)}) becomes
\begin{eqnarray}
E[X_t| \mathcal{F}_t^N]
&=&\frac{1}{1+{q_t}\dot{Q}_t}I(N_t=0),\label{(21a)}
\end{eqnarray}
%
where $q_t$ is defined as before with $\mu_0=0$, and
\begin{equation}\label{(19aa)}
\dot{Q}_t= {\int_{\{D\in\mathcal{U}:E_t\subseteq D\}}}e^ { \mu_1 |A_t\setminus D|}\,d\mu_\xi(D).
\end{equation}
It is easy now to see that the statement of Theorem~\ref{main} is
still valid in this case,
with $V$ replaced by $\dot{V}$, where
\[
\dot{V}_t=-c_1+(c_0+c_1+k_1\lambda_t)\frac{1}{1+q_t\dot{Q}_t}I(N_t=0).
\]
\end{Com}

\section{Examples}\label{Examples}

In this section, we apply Theorem~\ref{main} to our two
examples. We will see that in some sense they are are both analogous to the
univariate model of \cite{HJ}, in which the change-point is
exponentially distributed. There are two natural generalizations in
$\mathbf{R}_+^2$: first, $L$ is the single jump process in which the components
of the jump are independent univariate exponential random variables,
and second, $L$ is the first line of a Poisson process, noting that an
exponential
random variable can be regarded as the ``first line'' of a Poisson
process on $\mathbf{R}_+$. Although at first glance the single jump
process looks more straightforward, we shall see that in fact the
analysis is far more complex than in the case of the first line of a
Poisson process.
\begin{Ex}[(The single jump process)]\label{singlejump2}
Referring to Example~\ref{singlejump}, we have
$\lambda_t=\frac{f_t}{1-F_t}$ and
$q_t=
\frac{e^{-(\mu_1-\mu_0)t_1t_2}}{1-F_t}$.
Here we will consider the case in which the components $(Y_1,Y_2)$ of
the jump $Y$ are independent identically distributed exponential
random variables with parameter $\gamma$.
In this case,
\begin{eqnarray*}
\lambda_t&=&\frac{f_t}{1-F_t}=\frac{\gamma e^{-\gamma t_1}
\gamma e^{-\gamma t_2}}{1-(1-e^{-\gamma t_1})(1-e^{-\gamma t_2})}
=\frac{\gamma^2}{e^{ \gamma t_1} +e^{ \gamma t_2}-1},
\end{eqnarray*}
and is decreasing in each component. Next, we consider $q_t$:
\begin{eqnarray*}
q_t&=&\frac{e^{-(\mu_1-\mu_0)t_1t_2}}
{1-(1-e^{-\gamma
t_1})(1-e^{-\gamma t_2})}=\frac{e^{-(\mu_1-\mu_0)t_1t_2}} {
e^{-\gamma t_1}+e^{-\gamma t_2}- e^{- \gamma(t_1 + t_2 )}}.
\end{eqnarray*}
To find sufficient conditions to ensure that $q$ is increasing in
$t_1$ and $t_2$ on some set
$R=[0,r]^2$, we will assume that $\gamma>\mu_1-\mu_0$ and to
simplify the discussion (without loss of generality, by suitably
rescaling the time parameters if necessary) that
$\mu_1-\mu_0=1$. Now rewrite $q_t=1/g_t$
where
\[
g_t=g_{(t_1,t_2)}=e^{-t_1(\gamma-t_2)}(1-e^{-\gamma t_2})+e^{
t_1t_2}e^{-\gamma t_2}.
\]
We will show that if $r\leq\frac{\ln
\gamma}{\gamma} $, then $\frac{d}{dt_1}g_{(t_1,t_2)}\leq0$ for
$(t_1,t_2)\in R=[0,r]^2$. By symmetry, the same is true for
$\frac{d}{dt_2}g_{(t_1,t_2)}$ for $t\in R$. Therefore, $g$ is
decreasing and
$q=1/g$ is increasing in each component on $R$, and an optimal
solution exists for the sequential detection model.

To complete the example, we observe that
\begin{eqnarray*}
\frac{d}{dt_1}g_{(t_1,t_2)}= e^{-t_1(\gamma-t_2)}
\bigl(-(\gamma-t_2)(1-e^{-\gamma t_2})+t_2e^{\gamma(t_1-t_2)}\bigr)
\leq0
\end{eqnarray*}
if and only if
\[
(\gamma-t_2)(1-e^{-\gamma t_2})\geq t_2e^{\gamma(t_1-t_2)}
\]
or
equivalently,
%
\begin{equation}\label{(25)}
e^{\gamma t_1}\leq\frac{\gamma-t_2}{t_2}(e^{\gamma t_2}-1).
\end{equation}
The left-hand\vspace*{-1pt} side of (\ref{(25)}) is bounded above by $\gamma$
since $t_1\leq r\leq\frac{\ln\gamma}{\gamma}$. The right-hand side\vspace*{-1pt}
of (\ref{(25)}) is bounded below by
$\gamma$ since $t_2\leq r\leq\frac{\ln\gamma}{\gamma}\leq
\gamma-1$ when
$\gamma\geq1$, and so
$\frac{\gamma-t_2}{t_2}(e^{\gamma
t_2}-1)=\frac{e^{\gamma
t_2}-1}{t_2}(\gamma-t_2)\geq\gamma(\gamma-t_2)\geq
\gamma$. Therefore,\vspace*{-1pt} it is sufficient that $t_1,t_2\leq r\leq\frac
{\ln\gamma}{\gamma}$.
\end{Ex}

\begin{Ex}[(First line of a Poisson process)]\label{poisson2}
From the discussion in Example~\ref{poisson}, if $L$ is the
first line of a Poisson process with rate $\gamma$, then
$\lambda\equiv\gamma$, and so trivially is decreasing in each
component. We have
$q_t=
\frac{e^{-(\mu_1-\mu_0)t_1t_2}}{e^{-\gamma
t_1t_2}}=e^{(\gamma-(\mu_1-\mu_0))t_1t_2}$, which is increasing in each
component if $\gamma\geq\mu_1-\mu_0$. Therefore, an optimal
solution to
the sequential detection problem exists on any bounded set $R=[0,r]^2$
if
$\gamma\geq\mu_1-\mu_0$, and is defined by (\ref{(11)}). In fact, this
is exactly the same as the sufficient condition for the univariate
detection problem proven in
\cite{HJ} and~\cite{JH}.
\end{Ex}

\section{Conclusion}\label{Conclusion}

As indicated in the Introduction, the sequential detection model
considered here is only one of many scenarios that should be analyzed
in the general context of the ``optimal set-detection problem.''
Indeed, the model can be extended in many possible ways.
\begin{itemize}
\item The information structure: In addition to the sequential
information model, Herberts and Jensen
\cite{HJ} consider what they call the ``ex-post'' analysis. This would
correspond to observing $N$ on all of $R$, and then trying to
optimize the expectation of the valuation function. (Formally, this
corresponds to $\mathcal{F}_t=\mathcal{F}_{(r,r)} \,\forall t\in R$.)
Several variants
or combinations of the ex-post and sequential schemes can be studied.
\item The underlying space: We worked here
on a bounded subset of $\ R_{+}^{2} $. It would be of interest to
consider change-point problems on higher-dimensional Euclidean spaces
or more general partially ordered sets as in \cite{IMP}.
\item The change mechanism: Here the change occurs at either a
single random point or at the first line of a more general point
process. The example involving the first line of a Poisson process
turned out to be (perhaps surprisingly) the more natural analog of
the one-dimensional exponential change-point problem. Consideration
should be given to more general single jump and first line processes,
as well as to more general random sets (not necessarily upper layers).
For example, the case in which $L$ is the first line of an
inhomogeneous Poisson process with intensity $\gamma(\cdot)$ is
considered in \cite{Col} where it is proven that an optimal solution
exists if $\inf_{u\in R}\gamma(u)\geq\mu_1-\mu_0$.
\item The observed process: On $\mathbf{R}_+$, the process subject to
the change
can be a more general process, such as
the Brownian motion process (cf.~\cite{CMS}). Here too, we can
consider more general processes such as the set-indexed Brownian
motion (cf.~\cite{IM1}).
\item The parameters: In our analysis, it is implicitly assumed that
the parameters of the various processes are all known. How does one
approach the problem when one or more parameters must be estimated?
\item The gain function: Different valuation functions can be
chosen, thereby changing the notion of optimality. For example, with
a change generated by a single jump at $Y$, instead of two
cost parameters $c_0$ and $c_1$ associated respectively with $E_Y^c$
and $E_Y$, we could have different costs in each of the four
quadrants defined by
$Y$. Another variation considered in \cite{Col} is to replace $L_t$ in
(\ref{(5)}) with $I(L_t>0)$. Although this does not change the
valuation when the change is generated by a single jump, the analysis
becomes more complex when $L$ is the first line of a Poisson process.
\item Number of changes: Here we deal with only one change-set.
However, we can imagine that several changes occur on a decreasing
sequence of random upper layers, for example. This would correspond
to multiple change points on~$\mathbf{R}_+$.\looseness=1
\end{itemize}

\section*{Acknowledgments}

The second author thanks G. Ivanoff for her optimal hospitality while
visiting the University of Ottawa.

%

\printaddresses

\end{document}